\documentclass[11pt]{article}

\usepackage{amssymb,amsmath}

\usepackage[francais,english]{babel}
\usepackage[T1]{fontenc}

 \makeatletter \@addtoreset{equation}{section}
	
	\usepackage{hhline} 
	
	\usepackage{graphicx}
	
	\usepackage[english]{babel}
	\usepackage[english]{babel}
	\usepackage{amsfonts}
	\usepackage{xcolor}
	\usepackage{marginnote}
	\usepackage{amsmath}
	\usepackage{a4wide}
	\numberwithin{equation}{section}
	
	\usepackage{accents}

	\usepackage{amssymb}
	\usepackage{amsthm}
	
	\newtheorem{thm}{Theorem}[section]
	\newtheorem{lem}{Lemma}[section]
	\newtheorem{cor}{Corollary}[section]
	\newtheorem{prop}{Proposition}[section]
	
	\newtheorem{rem}{Remark}[section]

\begin{document}

		\title{Existence and uniqueness of a saddle-node bifurcation point for nonlinear equations in general domains}
		

		\author{Yavdat Il'yasov\\ \\
			Institute of Mathematics of UFRC RAS,   Ufa,
			Russia\\ ilyasov02@gmail.com}
		
		\date{$\ $}
		\maketitle

\begin{abstract}
This paper provides a direct method of establishing the existence and uniqueness of saddle-node bifurcations for  nonlinear equations in general domains. The method employs the scaled extended quotient whose saddle points correspond to the saddle-node bifurcations. The uniqueness of the saddle-node bifurcation point directly stems from the uniqueness of the saddle point. The method is applied to solving  open problems involving the existence and uniqueness of the maximum saddle-node bifurcation for the positive solutions curve  to an elliptic boundary value problem with a convex-concave nonlinearity in general domains.

\end{abstract}

\section{Introduction}
We study the saddle-node bifurcation points for the following problem
\begin{equation}
\label{f} 
\left\{ \begin{aligned}
		  -&\Delta u= \lambda u^{q}+ u^{\gamma},~~ x \in \Omega, \\[0.4em]
	~	&u|_{\partial \Omega}=0. 
\end{aligned}\right.
	\end{equation}
Here  $\Omega $ is a bounded connected domain in $\mathbb{R}^N$, $\lambda$ is a real parameter, $0<q<1<\gamma<2^*-1$ and $q\leq 2/(N-2)$ if $N\geq 3$. Here $2^*=2N/(N-2)$ if $N\geq 3$, $2^*=+\infty$ if $N=1,2$. We use the notations $W^{1,2}_0:= W^{1,2}_0(\Omega)$,  $F(u,\lambda):= -\Delta u- \lambda  u^{q}- u^{\gamma}$, $u \in W^{1,2}_0$, $\lambda \in \mathbb{R}$. 
A point $(u, \lambda) \in W^{1,2}_0\times \mathbb{R}$ is called a weak solution (or solution for short) of \eqref{f} if the equality $F(u, \lambda)=0$ is satisfied in the dual space $(W^{1,2}_0)^*$ of $W^{1,2}_0$. 
Define
$
S:=\{u \in C^2(\Omega)\cap C^1(\overline{\Omega}) \mid -\Delta u >  0,~u>  0 ~\mbox{in}~\Omega, ~u|_{\partial\Omega}=0\}$. 
Below we show that the map $F(\cdot,\lambda): W^{1,2}_0 \to (W^{1,2}_0)^*$ is continuously Fr\'echet differentiable on $S$, $\forall \lambda \in \mathbb{R}$.

We call a solution $(\hat{u},\hat{\lambda}) \in S\times \mathbb{R}$ of \eqref{f} the \textit{saddle-node  bifurcation point} (turning point, fold bifurcation) in $S$    if the following is fulfilled: (i) the nullspace $N(F_u(\hat{u},\hat{\lambda}))$ of the Fr\'echet derivative $F_u(\hat{u},\hat{\lambda})$  is not  empty; (ii) there exists $\varepsilon_1 >0$  and a neighbourhood $U_1 \subset W^{1,2}_0$  of $\hat{u}$ such that   for any $\lambda \in (\hat{\lambda},\hat{\lambda}+\varepsilon_1)$ equation \eqref{f} has no solutions in $S\cap U_1$; (iii) there exists $\varepsilon_2 >0$  and a neighbourhood $U_2 \subset W^{1,2}_0$  of $\hat{u}$ such that for each $\lambda \in (\hat{\lambda}-\varepsilon_2, \hat{\lambda})$   the equation has precisely two distinct solutions in $S\cap U$, cf. \cite{keller1, Kuznet}. This definition corresponds to the solutions curve turned to the left from the bifurcation value $\hat{\lambda}$. The solutions curve turned to the right  is defined similar.  In the case only (i)-(ii) are satisfied we call $(\hat{u},\hat{\lambda})$ the \textit{saddle-node type bifurcation point} of \eqref{f} in $S$. 
 A saddle-node  bifurcation point   $(u^*,\lambda^*)$ is said to be \textit{maximal} in $S$ if  $\hat{\lambda}\leq  \lambda^*$ for any other saddle-node  bifurcation $(\hat{u},\hat{\lambda})$ of \eqref{f} in $S$. We call a solution $u_\lambda \in C^2(\Omega)\cap C^1(\overline{\Omega})$ of \eqref{f}  \textit{stable} if $\lambda_1(F_u(u_\lambda,\lambda))\geq 0$, and \textit{asymptotically stable} if $\lambda_1(F_u(u_\lambda,\lambda))>0$, cf. \cite{Cazenave, CranRibin}. Here $\lambda_1(F_u(u_\lambda,\lambda))$ denotes the first eigenvalue of $F_u(u_\lambda,\lambda)(\cdot)$.
 
 
 

 
 The most widely used approach for detecting bifurcation points is the continuation of solutions by  a parameter (see, e.g., \cite{kielh,Korman2}). This method involves finding a limit point on the curve of solutions and testing if the limit point is a saddle-node type bifurcation. The failure of applicability of the Implicit Function Theorem is often used as a testing criterion.  In order to use this method, it is typically assumed that the functional $F(u, \lambda)$ at the bifurcation point is twofold differentiable.
 

 
 
 
 In this paper, we develop an alternative method introduced  in \cite{IlyasJDE24}, which consists of directly searching for the saddle-node bifurcation point. The key role in this method is played by the scaled extended quotient $\lambda(u,v)$, the saddle points of which correspond to the saddle-node bifurcations of the equation. Moreover, the uniqueness of a saddle point of its leads to the uniqueness of the saddle-node bifurcation point of the equation. An additional feature of this approach is that the conditions for smoothness on the functional $F(u, \lambda)$ are less restrictive than in the continuation methods.



First fundamental results on the structure of the set of positive solutions to \eqref{f} were obtained by A. Ambrosetti, G. Bresis, and H. Cerami in  \cite{ABC} (1994).
In particular, under the assumption $0<q<1<\gamma<+\infty$, they showed  that there exists $\Lambda>0$ such that (i) for all $\lambda \in (0, \Lambda)$ equation \eqref{f} has a stable minimal positive solution $u_\lambda$; (ii) for $\lambda=\Lambda$, equation \eqref{f} has at least one  positive solution $\hat{u} \in W^{1,2}_0\cap L^{\gamma+1}$; (iii) for all $\lambda>\Lambda$ equation \eqref{f} has no positive solutions. In the case $0<q<1<\gamma<2^*-1$, they showed  that for all $\lambda \in (0, \Lambda)$ there exists a second positive solution $w_\lambda$ of equation \eqref{f} such that $w_\lambda>u_\lambda$.

In the case $\Omega=B$  is the ball in $\mathbb{R}^n$, $0<q<1<\gamma<2^*-1$ the existence of the saddle-node  bifurcation point and the exact shape of the positive solutions curves of \eqref{f} was obtained  by P. Korman \cite{Korman} (2002),  M. Tang \cite{Tang} (2003), and  T. Ouyang, J. Shi \cite{Ouyang} (1999) under the assumptions $0<q<1<\gamma<\frac{n}{n-2}$. According to this result, equation \eqref{f} has exactly one positive solution for $\lambda=\lambda^*$, and exactly two positive solutions for $\lambda\in (0,\lambda^*)$. Moreover, all positive solutions lie on a single smooth solution curve, which for $\lambda<\lambda^*$ has two branches.  Much of the advancement in this result is largely because according to the famous result B. Gidas, W.-M. Ni, and L. Nirenberg theorem in \cite{Gidas}, any positive solution to elliptic equations on a ball must exhibit radial symmetry, which enables one to utilize the extensive tools available in ODE's theory.

The situation for general domain $\Omega$ is more
involved. The absence of results on the uniqueness of general-type solutions to \eqref{f} and the smoothness violation of the function $u^q$ at zero are major obstacles to such problems. This makes it challenging to justify the uniqueness of the bifurcation point and verify that solution curve \eqref{f} turns back from point $\lambda^*$.


Our approach relies on use of the \textit{scaled extended quotient} associated with equation \eqref{f}
\begin{equation}\label{eq:scaledQ}
	\lambda(u,v):=c_{\gamma,q}\frac{\langle-\Delta u, v\rangle^\frac{\gamma-q}{\gamma-1}}{\left\langle u^q,v\right\rangle\left\langle u^\gamma,v\right\rangle^\frac{1-q}{\gamma-1}},~~~~(u, v) \in S\times S
\end{equation}
where
$c_{\gamma,q}=\left(\frac{(\gamma-1)}{(\gamma-q)}
\frac{(1-q)}{(\gamma-q)} \right)^{(1-q)/(\gamma-1)}$.  We say that \textit{the minimax principle} is satisfied for $\lambda(u,v)$ in $S\times S$  if 
\begin{align}\label{Mmin}
	-\infty<	\lambda^*:=\sup_{ u \in  S}\inf_{v \in S}\lambda(u,v)=\inf_{v \in S}\sup_{ u \in  S}\lambda(u,v)<+\infty. 
\end{align} 
A pair  $(u^*,v^*) \in \overline{S}\times \overline{S}$ is called the \textit{saddle point} at $\lambda^*$  of $\lambda(u,v)$ in $S\times S$ if
\begin{equation}\label{sedsadd}
	\lambda^*=\lambda(u^*,v^*)=\sup_{ u \in  S}\lambda(u,v^*)=\inf_{v \in S}\lambda(u^*,v).
\end{equation}

Our main result is
\begin{thm}\label{thm1} Assume that $0<q<1<\gamma<2^*-1$ and $q\leq 2/(N-2)$ if $N\geq 3$. Then 
  \eqref{f} has at most one maximal saddle-node bifurcation point $(u^*,\lambda^*)$ in $S$.
	Moreover, $0<\lambda^*<+\infty$, and
	\begin{description}
	\item[$(1^o)$] $u^*$ is a stable positive solution of \eqref{f}, and $u^* \in S$, moreover there exists $v^* \in S$ such that  $N(F_u(\hat{u},\hat{\lambda}))=\{v^*\}$;
	\item[$(2^o)$]	the minimax principle
	is satisfied for $\lambda(u,v)$ in $S\times S$ and $(u^*,v^*) \in S\times S$ is a \textit{saddle point}  at $\lambda^*$  of $\lambda(u,v)$ in $S\times S$;
	\item[$(3^o)$] 
	the solution curve $(u_\lambda, \lambda)$ continues smoothly through $(u^*,\lambda^*)$ in such a way that for each $\lambda \in (\lambda^*-\varepsilon,\lambda^*)$ with some $\varepsilon>0$ equation \eqref{f} has two distinct positive solutions in $S$; 
	\item[$(4^o)$] for any $\lambda>\lambda^*$ equation \eqref{f} has no non-negative solutions;
	\item[$(5^o)$] for any $\lambda \in (0,\lambda^*)$ equation \eqref{f} possesses a  stable positive solution $u_\lambda \in S$.
\end{description}
	
\end{thm}

 \begin{rem}
 	In numerical methods, continuation via a parameter is commonly used to find bifurcations \cite{keller1, Kuznet}.  Minimax formulas of the extended quotient like \eqref{Mmin}, being variational, make it possible to develop numerical algorithms based on other principles for finding bifurcations, such as gradient-based algorithms or others,  see, e.g.,  \cite{Bagirov, IvanIlya, IlIvan1}. These algorithms may have certain advantages (see, e.g., \cite{ IvanIlya, IlIvan1, Salazar}) over the commonly used ones  \cite{keller1, Kuznet}.
 	The obtained minimax formula \eqref{Mmin} is expected to enhance the ability to numerical find bifurcations.
 	
 \end{rem}
 
\begin{rem}
 From \cite{IlyasJDE24} it follows that there exists $\bar{\lambda}\leq \lambda^*$ such that $\forall \lambda \in (0,\bar{\lambda})$ equation \eqref{f} has an asymptotically stable  positive solution $u_\lambda \in   S$.
\end{rem}

\section{Preliminaries}
If $u \in \overline{S}$, then the maximum principals for the elliptic problems imply that  $u>0$ in $\Omega$, and $\displaystyle{\min_{x' \in \partial \Omega}|\frac{\partial u(x')}{\partial \nu(x')|} >0}$, where $\nu(x')$ denotes the interior unit normal at $x' \in \partial \Omega$, see, e.g., Theorem 3.1, Lemma 3.4 in \cite{Giltrud}. 
Thus,
$$\displaystyle{\overline{S} \subset \tilde{S}:= \{u \in C^2(\Omega)\cap C^1(\overline{\Omega})\mid~u > 0 ~\mbox{in}~\Omega,~\min_{x' \in \partial \Omega}|\frac{\partial u(x')}{\partial \nu(x')|} >0,  ~u(x')=0,~x' \in \partial \Omega\}}.
$$
Define  $
S(\delta):=\{u \in  S \mid~~\min_{x' \in \partial \Omega}|\partial u(x')/\partial \nu(x')| >\delta\}$,  $S(0)=S$ for $\delta\geq 0$.  
Let $Y$ be a topological space such that $S(\delta) \subset Y$.  The set $S(\delta)$ endowed with topology of $Y$ we denote by $S_{Y}(\delta)$. $\mathcal{L}(S_{Y}(\delta),(W^{1,2}_0)^*)$ denotes the Banach space of bounded linear operators from $S_{Y}(\delta)$ into  $(W^{1,2}_0)^*$. 
  
  \begin{prop}\label{prop:cont} Let $0<q\leq 2/(N-2)$, $q<1<\gamma<2^*-1$, $\lambda \in \mathbb{R}$. 
  	The map $F(\cdot,\lambda):S_{W^{1,2}_0} \to (W^{1,2}_0)^*$ is Fr\'echet differentiable, moreover, $F_u(u,\lambda) \in C(S_{C^1}; \mathcal{L}(S_{C^1},(W^{1,2}_0)^*))$ and $F_u(u,\lambda) \in C(S_{L^p}(\delta); \mathcal{L}(S_{L^p}(\delta),(W^{1,2}_0)^*))$ for any $\delta>0$ and $p=Nq $. 
  	\end{prop}
  \begin{proof}  We develop the approach proposed for the proof of Lemma 3.2 in \cite{ABC}. Obviously, it suffices to verify that the assertions of the proposition are satisfied for the map $Q(u):=u^q$, $u \in W^{1,2}_0$.
  Since	$S \subset \tilde{S}$, there exists $c(u) >0$,  which does not depend on $ x \in \overline{\Omega}$ such that $u(x)\geq c(u)\mbox{ dist}(x,\partial \Omega)$ in $\Omega$.
  	By H\"older's and Sobolev's inequalities we have
  	$$
  	\int u^{q-1}\phi\psi  \, dx =\int u^q\left(\frac{\phi}{u}\right) \psi \, dx\leq \frac{1}{c(u)}\|u\|^q_{p}\cdot\|\frac{\phi}{\mbox{ dist}(\cdot,\partial \Omega)}\|_2\cdot \|\psi\|_{2^*}, ~~~\forall \phi, \psi\in W^{1,2}_0,
  	$$
  	where $p=Nq$. Note that $p\leq 2^*$ since $q\leq 2/(N-2)$.    	By the Hardy inequality $\|\phi/\mbox{ dist}(\cdot,\partial \Omega)\|_2\leq C\|\phi\|_{1,2}$,  $\forall \phi \in W^{1,2}_0$, and therefore, using Sobolev's inequalities we derive
  	\begin{equation}\label{welldef}
  		\int u^{q-1}\phi\psi \, dx \leq \frac{C}{c(u)} \|u\|^q_{p
  		}\|\phi\|_{1,2}\|\psi\|_{1,2}\leq \frac{C'}{c(u)} \|u\|^q_{1,2
  		}\|\phi\|_{1,2}\|\psi\|_{1,2}, ~~~\forall \phi,\psi\in W^{1,2}_0,
  	\end{equation}
  	for some $C,C'<+\infty$ which does not depend on $u, \phi, \psi\in W^{1,2}_0$. Hence, $Q(\cdot):S_{W^{1,2}_0} \to (W^{1,2}_0)^*$ is Fr\'echet differentiable.

  	It is easily seen that $c(\cdot) \in C(S_{C^1})$ and $\inf\{c(u)\mid~u \in S(\delta)\}=: c_1(\delta)>0$.   	
  	This by \eqref{welldef} yields that $Q_u(\cdot) \in  C(S_{C^1}; \mathcal{L}(S_{C^1},(W^{1,2}_0)^*))$ and $Q_u(\cdot) \in C(S_{L^p}(\delta); \mathcal{L}(S_{L^p}(\delta),(W^{1,2}_0)^*))$ for any $\delta>0$ and  $p=Nq$.    
  \end{proof}

Observe that the saddle-node bifurcation point $(u, \lambda) \in S\times \mathbb{R}$ of \eqref{f} should to satisfy the system of equations
\begin{equation}\label{eq:Msys} 
		\left\{ \begin{aligned} 	&F(u, \lambda)   =0,\\ 	&F_u(u, \lambda) (v)=0, 	
	\end{aligned}\right.
\end{equation}
with some $v \in N(F_u(u, v) )$.
The extended Rayleigh quotient  \cite{IlyasJDE24} associated with \eqref{f} is defined as follows
\begin{equation}
	R(u,v):=\frac{\left\langle -\Delta u, v\right\rangle-\left\langle u^{\gamma},v\right\rangle}{\left\langle u^q,v\right\rangle}, ~~ (u, v) \in S\times S.
\end{equation}
Observe,  
\begin{align}
	\left\{ \begin{aligned}
		&\lambda=R(u,v),\\
		&R_v(u,v)=0,\\
		&R_u(u,v)=0
			\end{aligned}\right.~\Leftrightarrow~
			\left\{ \begin{aligned} 	&F(u, \lambda)   =0,\\ 	&F_u(u, \lambda) (v)=0, 	
			\end{aligned}\right.
\end{align}
Thus the set of solutions of \eqref{eq:Msys} coincides with the set of critical points of $R(u,v)$ on $S\times S$.

Consider
$$
R(tu,v)=\frac{t^{1-q}\left\langle -\Delta u, v\right\rangle-t^{\gamma-q}\left\langle u^{\gamma},v\right\rangle}{\left\langle u^q,v\right\rangle},~~t>0,~ (u, v) \in S\times S.
$$
Note that $\left\langle u^q,v\right\rangle \neq 0$, $\left\langle u^\gamma,v\right\rangle  \neq 0$, $\forall  (u,v) \in S \times S$ since $u,v>0$ in $\Omega$.
It is easily seen that for any $u,v \in S$, the function $R(tu,v)$ defined on  $t\geq 0$  has a  absolute maximum at the  point 
\begin{equation}\label{t}
	t(u,v)=\left(\frac{(1-q)\left\langle -\Delta u, v\right\rangle}{(\gamma-q)\left\langle u^{\gamma},v\right\rangle}\right)^{1/(\gamma-1)}>0.
\end{equation}
 It is easily to check that the scaled quotient \eqref{eq:scaledQ} is defined as 
$$
\lambda(u,v)=R(t(u,v)u,v)=\max_{t\geq 0}R(tu,v),~~(u, v) \in S\times S.
$$
Observe, $R_t(tu,v)|_{t=t(u,v)}=0$. This implies that  $\lambda_u(u,v)=0~\Leftrightarrow~R_u(u,v)=0 $, and 
$\lambda_v(u,v)=0~\Leftrightarrow~R_v(u,v)=0$.

\begin{prop}\label{prop1}
If 	$u \in  W^{1,2}_0$ is a weak non-negative solution   to \eqref{f}, then $u \in  S$. 
\end{prop} 
\begin{proof} Note that equality \eqref{f} for $u$ implies that $-\Delta u\geq 0$. 
	The standard bootstrap argument and  Sobolev's embedding theorem entail that  $u \in L^\infty(\Omega)$.   Therefore, by the $L^p$-regularity results from \cite{Giltrud}, $u \in W^{2,p}(\Omega)$ for any $1 < p < \infty$ and thus, by Sobolev's embedding theorem, $u \in C^{1,\alpha}(\overline{\Omega})$ for any $\alpha \in (0, 1)$. Moreover, Schauder estimates  (see, e.g., \cite{Giltrud}) imply that $u \in C^{2}(\Omega)$, and thus $u \in  S$.
	
\end{proof}

\section{Minimax principal}
A function 
$f:X\to \mathbb{R}$  defined on a convex subset 
$X$ of a real vector space is  \textit{quasi-convex} (strong quasi-convex) if for all 
$x,y\in X$ and 
$\lambda \in [0,1]$ we have
$$
f(\lambda x+(1-\lambda )y)\leq (<) \max {\big \{}f(x),f(y){\big \}}.
$$
A function whose negative is quasi-convex (strong quasi-convex) is called \textit{quasi-concave} (strong quasi-concave).
We need on  Sion's minimax theorem \cite{Sion}

\noindent
{\bf Theorem} (\textit{Sion})
\textit{Let $X$ be a compact convex subset of a linear topological space and 
	$Y$ a convex subset of a linear topological space. If 
	$f$ is a real-valued function on 
	$X\times Y$  with}
\begin{itemize}
	\item \textit{$f(x,\cdot)$ upper semicontinuous and quasi-convex on 
		$Y$, $\forall x\in X$, and}
	\item $f(\cdot,y)$\textit{ lower semicontinuous and quasi-concave on }
	$X$, $\forall y\in Y$,
\end{itemize}
\textit{then the minimax principle is satisfied}
$$
\min_{x\in X}\max_{y\in Y} f(x,y)=\max_{y\in Y}\min_{x\in X}f(x,y).
$$

\begin{lem}\label{lem3.1}
$\lambda(u,v)$ satisfies the minimax principle \eqref{Mmin} in $S\times S$. 
\end{lem}
\begin{proof}
The following lemma is proven below in  Appendix A.
\begin{lem}\label{propSion}
	\begin{description}
		\item[\rm{1)}] $\lambda(u,\cdot)$ is  quasi-convex function on $S$,
		$\forall u \in S$,
		\item[\rm{2)}] $\lambda(\cdot, v)$ is strong quasi-concave function on $S$, 
		$\forall v \in S$.
	\end{description}
\end{lem}
Denote by $\phi_1 \in S$  the eigenfunction corresponding to the principal eigenvalue $\lambda_1:=\lambda_1(-\Delta)$. Observer,
$$
\lambda^*:=\sup_{ u \in  S}\inf_{v \in S}\lambda(u,v)\leq \sup_{u \in S}\lambda(u, \phi_1)=\sup_{u \in S}\frac{\lambda_1 \int \left(u -u^{\gamma} \right) \phi_1}{\int u^{q}\phi_1}.
$$
Note that  
$\lambda_1u -u^{\gamma} \leq K u^{q}$, $\forall u \in \mathbb{R}$ for sufficiently large $K$. Hence, $\lambda^*\leq K<+\infty$.
From \cite{ABC} it follows that there exists $\Lambda>0$ such that \eqref{f} has a positive solution $u_\lambda\in S$ for any $\lambda \in (0,\Lambda)$. Since $\lambda=\inf_{v \in S}\lambda(u_\lambda,v)$, $0<\lambda\leq \lambda^*$, $\forall \lambda \in (0,\Lambda)$. We thus have $0<\lambda^*<+\infty$.

Let $\varepsilon >0$, $r>\varepsilon ~\|\phi_1\|_{C^1(\Omega)}$. Introduce 
$$
S(\varepsilon, r):=\{w \in S:~ \varepsilon \phi_1 \leq w~\mbox{in}~\Omega, ~\|w\|_{C^1(\Omega)}\leq r\}.
$$ 
Clearly, $\lambda(u,v)$  is a continuous functional on $S\times S$, and  $S(\varepsilon, r)$ is a compact convex subset in $C(\overline{\Omega})$. Hence, Sion's Theorem and Lemma \ref{propSion} yield 
\begin{equation}\label{eq:3.1}
-\infty<	\lambda^*(\varepsilon, r):=\sup_{ u \in  S(\varepsilon, r)}\inf_{v \in S}\lambda(u,v)=\inf_{v \in S}\sup_{ u \in  S(\varepsilon, r)}\lambda(u,v)\leq \lambda_1^*<+\infty,
\end{equation} 
for sufficiently small $\varepsilon$ and large $r$. Note that the positive solution $u_\lambda \in  S(\varepsilon, r)$ for any $\lambda \in (0,\Lambda)$ and sufficiently small $\varepsilon$ and large $r$, and therefore $0<\lambda\leq \lambda^*(\varepsilon, r)$, $\forall \lambda \in (0,\Lambda)$. 
Passing to the limit in \eqref{eq:3.1} as $\varepsilon \to 0$ and $r\to +\infty$ we obtain  \eqref{Mmin}. 
\end{proof}

\section{Existence of the saddle-node type bifurcation point }

Define $	S_s:=\{u \in S:~~\lambda_1(F_u(u,\lambda))> 0\}$, $\Sigma(u):=\{v \in W^{1,2}_0: \int uv   \neq 0\}$,  $u \in S$. Consider 
\begin{equation}\label{MainB}
	\lambda^*_s:= \sup_{u\in S_s}\inf_{v\in \Sigma(u)}R(u, v).
\end{equation}
Since $S\subset \Sigma(u)$ and $S_s \subset S$, $\lambda^*_s\leq  \sup_{u\in S_s}\inf_{v\in S}R(u, v)\leq \lambda^*<+\infty$. 
From \cite{IlyasJDE24} it follows that there exists $\bar{\lambda}>0$ such that \eqref{f} has an asymptotically stable  positive solution $u_\lambda\in S$ for any $\lambda \in (0,\bar{\lambda})$. Hence, $0<\bar{\lambda}\leq \lambda^*_s$. Thus, $0<\lambda^*_s\leq \lambda^*<+\infty$.
	\begin{lem}\label{lem:nonex}
	For  $\lambda>\lambda^*_s$,  \eqref{f} has no stable weak non-negative solutions.
\end{lem}
\begin{proof}
	Suppose by contradiction that for $\lambda>\lambda^*_s$ there exists a stable weak non-negative solution $u_\lambda$ of \eqref{f}. Using Proposition \ref{prop1} we obtain $u_{\lambda}\in S_s$,  and consequently,   \eqref{MainB} implies $\inf_{v\in \Sigma(u_\lambda)}\mathcal{R}(u_\lambda, v)<\lambda$. Hence,  $\exists v \in \Sigma(u_\lambda)$ such that $\mathcal{R}(u_\lambda, v)<\lambda$. But this contradicts to the equality $F(u_\lambda, \lambda)=0$.
\end{proof}

\begin{prop}\label{lem1} Let $u \in S_s$. Then  
	$$
	-\infty<\lambda(u):=\inf_{v \in \Sigma(u)}R(u ,v)<+\infty
	$$
	if and only if 	$u$ is a   weak positive solution of \eqref{f} with $\lambda=\lambda(u)$.
	\end{prop}
	\begin{proof}
		 Assume that $u$ is a weak positive solution of \eqref{f} with $\lambda=\lambda(u)$. Then $R(u, v)=\lambda(u)$ for all $v\in \Sigma(u)$, and thus, we get $\lambda(u):=\inf_{v \in \Sigma(u)}R(u ,v)$. The reverse assertion follows from Lemma 3.2. in  \cite{IlyasJDE24}.
	\end{proof}

\begin{lem}\label{lemM} 
	 Equation \eqref{f}
	for $\lambda=\lambda^*_s$  has a stable positive solution $u^* \in   \overline{S}_s$. 

	\par 
	
\end{lem}

\begin{proof}
	
 Since $-\infty<\lambda^*_s<+\infty$,  there exists a sequence 
		$u_n \in S_s$, $n=1,\ldots$, such that 
		$$
		\lambda_n:=\lambda(u_n):=\inf_{v \in \Sigma(u_n)}\mathcal{R}(u_n,v) \to \lambda^*_s~~\mbox{as}~~n\to +\infty.
		$$ 
	By Proposition \ref{lem1}  
	\begin{equation}\label{BEq1}
		-\Delta u_n = \lambda_n  u_n^q+ u_n^\gamma, ~~n=1,\ldots.
	\end{equation}
Since $\lambda_1(F_u(u_n,\lambda))> 0$, $\langle F_u(u_n,\lambda)(u_n),u_n \rangle>0$.  This easily implies (see \cite{IlyasJDE24}) that  sequence $\|u_n\|_{1,2}$ is bounded. 
	Hence by the Banach–Alaoglu theorem  there exists a subsequence (again denoted by $(u_n)$) such that
	\begin{equation}\label{wse}
		u_n \to u^*~~\mbox{weakly in}~  W^{1,2}_0 ~~\mbox{and strongly in}~~ L^\beta~ ~\mbox{as}~ n\to +\infty,
	\end{equation}
 for some $u^* \in  W^{1,2}_0$ and $1\leq \beta < 2^*$.
	To show that $ u^*\neq 0$, we  use the following lemma from \cite{ABC}
	\begin{prop} \label{ABCLem}
		Assume that $f(t)$ is a function such that $t^{-1}f(t)$ is decreasing for $t>0$. Let $v$ and $w$ satisfy:
		$u>0$, $w>0$ in $\Omega$, $v=w=0$ on $\partial \Omega$, and
		\begin{align*}
			-\Delta w\leq f(w),~~-\Delta u\geq f(u), ~~\mbox{in}~~ \Omega.
		\end{align*}
		Then $u\geq w$.
	\end{prop}
	Let $w_\lambda$ be a unique positive solution of 
	\begin{equation}
		\label{q}
		\left\{ \begin{aligned}
			-\Delta w_\lambda=&\lambda w_\lambda^{q-1}~~\mbox{in}~~ \Omega,\\
			w_\lambda |_{\partial \Omega}&= 0.
		\end{aligned}\right.
	\end{equation}
	Note that $w_{\lambda}\equiv (\lambda)^{1/(2-q)}w_1$.
	By \eqref{BEq1},  $ -\Delta u_n  \geq  \lambda_n  (u_n)^{q-1} ~~  in ~~ \Omega$, $ n=1\ldots$.
Hence, Proposition  \ref{ABCLem} yields
	\begin{equation}\label{ineqABC}
		u_n \geq w_{\lambda_n}\equiv (\lambda_n)^{1/(2-q)}w_1, \quad i=1,\ldots,m,
	\end{equation}
	and consequently, $u^*\geq w_{\lambda^*_s}>0$.
	Passing to the limit in \eqref{BEq1} as $n\to +\infty$ we obtain
	\begin{equation*}
		-\Delta u^*  = \lambda^*_s  (u^*)^q + (u^*)^\gamma,~~x \in \Omega.
	\end{equation*}
	Using Proposition \ref{prop1} we conclude that   $u^* \in  S$. It remains to show that  $u^* \in \overline{S}_s$
	Clearly, $w_{\lambda^*_s}$ satisfies the Hopf boundary maximum principle. As a result,  $w_{\lambda^*_s} \in S(\delta(w_{\lambda^*_s}))$ with some $\delta(w_{\lambda^*_s})>0$. 
	Inequalities \eqref{ineqABC} imply that $u_n \in S(\delta(w_{\lambda^*_s}))$, $n=1,\ldots $. In virtue of Proposition \ref{prop:cont}, $F_u(u,\lambda) \in C(S_{L^p}(\delta(w_{\lambda^*_s})); \mathcal{L}(S_{L^p}(\delta(w_{\lambda^*_s})),(W^{1,2}_0)^*))$, with $p=Nq \in (1,2^*)$. Hence, 
		\begin{align}\label{ConvDR}
		\langle F_u(u_n,\lambda_n)(\phi),\psi \rangle \to
		\langle F_u(u^*,\lambda^*_s)(\phi),\psi \rangle~ \mbox{as}~~ n\to +\infty&,\quad\forall \phi,\psi  \in W^{1,2}_0.
	\end{align}
 Thus, $\lambda_1(F_u(u_n,\lambda_n)) \to \lambda_1(F_u(u^*,\lambda^*_s))\geq 0$ as $n\to +\infty$, and therefore, $u^* \in \overline{S}_s$.
\end{proof}
\begin{lem}\label{lem:spEx}
$\lambda_1(F_u(u^*,\lambda^*_s))=0$, that is, $N(F_u(u^*,\lambda^*_s))\neq \emptyset$.
\end{lem}
\begin{proof}
	Suppose, contrary to our claim, that $\lambda_1(F_u(u^*,\lambda^*_s))>0$, that is,   $F(\cdot,\lambda^*_s):S_{C^1} \to (W^{1,2}_0)^*$  is nonsingular linear operator. In virtue of Proposition \ref{prop:cont}, $F(\cdot,\lambda) \in C^1(S_{C^1};\mathcal{L}(S_{C^1},(W^{1,2}_0)^*))$. 
	Hence, by the Implicit Functional Theorem (see, e.g, \cite{Dieudonn}) there is a neighbourhood $V_1\times U_1 \subset V\times U$  of $(\lambda^*_s, u^*)$ in $\mathbb{R}\times S_{C^1}$ and a mapping $V_1 \ni \lambda \mapsto u_\lambda \in U_1$ 
such that $u_\lambda|_{\lambda=\lambda^*_s}=u^*$ and $F(u_{\lambda}, \lambda)=0$, $\forall \lambda \in V_1$. 
Furthermore, the map $u_{(\cdot)}: V_1 \to X$ is continuous. Since $\lambda_1(F_u(u^*,\lambda^*_s))>0$, there exists a neighbourhood $V_2 \subset V_1$ of $\lambda^*_s$ such that $\lambda_1(F_u(u_{\lambda},\lambda))>0$ for every $\lambda \in V_2$. However, this contradicts Lemma \ref{lem:nonex}, and thus,  $\lambda_1(F_u(u^*,\lambda^*_s))=0$.  

\end{proof}
From Lemmas \ref{lem:nonex}, \ref{lemM},  \ref{lem:spEx} we obtain
\begin{cor}\label{cor2}
	$(u^*,\lambda^*_s)$ is a saddle-node type bifurcation point of \eqref{f} in $S$.
\end{cor}

\section{Proof of Theorems \ref{thm1}}

We claim that $\lambda^*_s=\lambda^*$. Suppose the converse $\lambda^*_s<\lambda^*$. Then for any $\lambda \in (\lambda^*_s,\lambda^*)$ there exists $\check{w}_\lambda \in S$ such that $R(\check{w}_\lambda, v)>\lambda$, and consequently $\langle F(\check{w}_\lambda,\lambda), v \rangle >0 $ for all $v \in S$, which means that $\check{w}_\lambda$ is a supersolution of \eqref{f}. It is easily seen that for sufficiently small $\varepsilon$, the function $\varepsilon \phi_1$ is a subsolution of \eqref{f} and $\varepsilon \phi_1 <\check{w}_\lambda$ in $\Omega$. Hence by the super-sub solution method  there exists 
a minimal solution $u_\lambda \in S$  of \eqref{f} such that $\varepsilon \phi_1\leq u_\lambda \leq \check{w}_\lambda$ (see, e.g., \cite{ABC}). Furthermore, Lemma 3.2 in \cite{ABC} implies that $\lambda_1(F_u(\check{w}_\lambda,\lambda))\geq 0$, and thus, $\check{w}_\lambda \in \overline{S}_s$. Notice 
$R(\check{w}_\lambda,v)=\lambda>\lambda_s^*$, $\forall v \in S$, and $\inf_{v\in S} R(u,v)<\lambda$, $\forall u \in S_c$, we obtain a contradiction. Hence, $\lambda^*_s=\lambda^*$.

Let us prove $(1^o)$.
By Corollary \ref{cor2}  there exists a saddle-node type bifurcation point $(u^*, \lambda^*) \in \overline{S}_s\times \mathbb{R}_+$ of \eqref{f}  such that 1)  $\forall \lambda>\lambda^*$ equation \eqref{f} has no stable positive   solutions, 2) $u^*$ is a stable positive solution of \eqref{f}, 3) $\lambda_1(F_u(u^*,\lambda^*))=0$, that is, $F_u(u^*,\lambda^*)$  is singular. This yields that there exists non negative $v^* \in  W^{1,2}_0\setminus 0$ such that
$$
-\Delta v^* -\lambda^* q (u^*)^{q-1}v^*-\gamma (u^*)^{\gamma-1}v^*= 0.
$$
Since $\lambda^* q (u^*)^{q-1}v^*+\gamma (u^*)^{\gamma-1}v^*\geq 0$ in $\Omega$, the maximum principle yields that $v^*>0$ in $\Omega$, and thus, $v^* \in S$. Consequently,  $N(F_u(\hat{u},\hat{\lambda}))=\{v^*\} \subset S$. Thus,  we get  $(1^o)$.

Lemma \ref{lem3.1} implies  
$\lambda^*=\lambda(u^*,v^*)=\sup_{ u \in  S}\inf_{v \in S}\lambda(u,v)=\inf_{v \in S}\sup_{ u \in  S}\lambda(u,v)$, that is, $(u^*,v^*) \in S\times S$ is a saddle point  of $\lambda(u,v)$ in $S\times S$, and thus, we get  $(2^o)$.

Since $\int (u^*)^q v^* >0$,  $F_\lambda(u^*,\lambda^*) = -(u^*)^q \notin R(F_u(\hat{u},\hat{\lambda}))$. Furthermore, in virtue of Proposition \ref{prop:cont} the map $F(\cdot,\cdot):S_{C^1}\times \mathbb{R}_+ \to (W^{1,2}_0)^*$  is  continuously differentiable in an open neighbourhood of $(u^*, \lambda^*)$. Hence, the Crandall \& Rabinowitz Theorem (see Theorem 3.2 in \cite{CranRibin} ) yields that  solutions of $F(u,\lambda)=0$ near $(u^*, \lambda^*)$ form a continuous  curve $(u(s),\lambda(s)) \in S\times \mathbb{R}_+$, $s \in (-a,a)$,  with some $a>0$ such that $(u(0),\lambda(0))=(u^*, \lambda^*)$ and $F(u(s),\lambda(s))=0$, $s \in (-a,a)$. The uniqueness of point $u^*$ implies that $\lambda(s)\neq\lambda^*$ for $\forall s \in (-a,a)\setminus 0$. Hence, there exists $\varepsilon>0$ such that for all $\lambda \in (\lambda^*-\varepsilon,\lambda^*)$ equation \eqref{f} has two distinct positive solutions in $ S$, and we obtain $(3^o)$.

Analysis similar to the proof of Proposition \ref{prop1} it follows that equation \eqref{f} has no nonnegative  solutions for $\lambda>\lambda^*$. Thus, assertion $(4^o)$ of the theorem holds true. 

By Theorem 2.1, Lemma 3.5 from \cite{ABC} it follows that there exists $\Lambda>0$ such that for any $\lambda \in (0,\Lambda)$ equation \eqref{f} possesses stable positive solution $u_\lambda \in C^2(\Omega)\cap C^1(\overline{\Omega})$, while for $\lambda >\Lambda$  equation \eqref{f} has no positive solutions. This and assertions $(1^o)$, $(3^o)$-$(4^o)$ of the theorem imply that  $\lambda^*=\Lambda$, and thus we obtain $(5^o)$.

Let us show that $u^*$ is uniquely defined. Suppose that there exists a second  saddle point  $(u^{**}, v^{**})$ at $\lambda^*$ of $\lambda(u,v)$ in $S\times S$, that is, $\lambda^*=\lambda(u^{**}, v^{**})=\sup_{u\in S_s}\lambda(u,v^{**})=\inf_{v \in S}\lambda(u^{**},v)$, and $u^*\neq u^{**}$. 
Observe,
\begin{equation*}
	\lambda^*=\inf_{v \in S}\lambda(u^*,v) \leq \lambda(u^{*}, v^{**})\leq \sup_{u\in S_s}\lambda(u,v^{**})=\lambda^*
\end{equation*}
Hence $\lambda(u^{*}, v^{**})=\lambda^*$. Since 
$\lambda(u,v)$ is strong quasi-convex with respect to $u$, 
$$
\lambda^*=\sup_{u\in S_s}\lambda(u,v^{*}) \geq\lambda(\tau u^*+(1-\tau)u^{**},v^*)> \min\{\lambda(u^*,v^{*}),\lambda(u^{**},v^{*})\}=\lambda^*,~~\tau \in [0,1].
$$
We get a contradiction, and thus, $u^*=u^{**}$. 

Since $u^*$ is uniquely defined, Corollary \ref{cor2} and assertions $(1^o)$-$(5^o)$ yield that $(u^*,\lambda^*)$ is a uniquely defined saddle-node  bifurcation point of \eqref{f} in $S$.

\section{Acknowledgment} The author is greatly indebted to Nurmukhamet Valeev for his collaboration in proving Lemma 3.2.

\section{Appendix A}

{\it Proof of Lemma \ref{propSion}:}
Denote $\beta_1:=\frac{\gamma-1}{\gamma-q}$, $\beta_2:=\frac{1-q}{\gamma-q}$. Let
$$
\tilde{\lambda}(u,v):=\frac{\left\langle -\Delta u, v\right\rangle}{\left\langle u^q,v\right\rangle^{\beta_1}\left\langle u^\gamma,v\right\rangle^{\beta_2}},~~u\in S, v \in S.
$$
Clearly, it is  sufficient to prove the following 
	\begin{description}
	\item[\rm{(i)}] $\tilde{\lambda}(\cdot, v)$ is a strong quasi-concave function on $S$, 
	$\forall v \in S$, 
	\item[\rm{(ii)}] $\tilde{\lambda}(u,\cdot)$ is a quasi-convex function on $S$, 
	$\forall u \in S$, 
\end{description}
Let us prove  (i). Consider
	\begin{align*}
		f(\alpha):=&\frac{\alpha\left\langle -\Delta  u, v\right\rangle+(1-\alpha)\left\langle -\Delta  w, v\right\rangle}{\left(\left\langle (\alpha u+(1-\alpha)w )^q,v\right\rangle\right)^{\beta_1}\left(\left\langle (\alpha u+(1-\alpha)w )^\gamma,v\right\rangle \right)^{\beta_2}}, ~~\alpha \in [0,1],~u, v \in S.
	\end{align*}
	To prove the assertion, we need to show
	$
		f(\alpha)> \min\{f(0),f(1)\}\equiv\min\{\tilde{\lambda}(u,v)),\tilde{\lambda}(w,v))\}
	$.
		Denote
	\begin{align*}
		a_1:=\left\langle -\Delta  u, v\right\rangle,~~ a_2:=\left\langle -\Delta  w, v\right\rangle,~~t=\alpha a_1+(1-\alpha)a_2.
	\end{align*}
Assume first that $a_1-a_2\neq 0$. Without loss of generality we may assume that $a_1-a_2> 0$.
	Observe
		\begin{equation}\label{talpha}
		\alpha u+(1-\alpha)w =t\frac{(u-w)}{a_1-a_2}+\frac{a_1w-a_2u}{a_1-a_2},
	\end{equation}
	and
	\begin{align*}
		\tilde{f}(t):=f(\frac{t-a_2}{a_1-a_2})=
		\frac{(a_1-a_2)}{\left(\left\langle \left((u-w)+\frac{a_1w-a_2u}{t} \right)^q,v\right\rangle\right)^{\beta_1}\left(\left\langle\left((u-w)+\frac{a_1w-a_2u}{t} \right)^\gamma,v\right\rangle \right)^{\beta_2}}
	\end{align*}
	Observe that if $a_1w-a_2u = 0$ in $\Omega$, then $\tilde{f}'(t)=0$, $\forall t \in [a_2,a_1]$, and thus the proof follows. 
	\begin{prop}\label{V1}
		Suppose that $a_1w(x)-a_2u(x) \not\equiv 0$  in $\Omega$.	If  $(a_1w-a_2u)(x)$  does not change  sign in $\Omega$, then $f(\alpha) > \min\{f(0),f(1)\}$, $\forall \alpha \in [0,1]$.
	\end{prop}
	\begin{proof}
		Observe
		\begin{align*}
			\frac{d}{dt}\psi_q(t):=&\frac{d}{dt}\left\langle\left((u-w)+\frac{a_1w-a_2u}{t} \right)^q,v\right\rangle^{\beta_1}=\\
			&-q\beta_1\frac{1}{t^3}\left\langle\left((u-w)+\frac{a_1w-a_2u}{t} \right)^{q-2}\left(t(u-w)+a_1w-a_2u \right)(a_1w-a_2u),v\right\rangle\times\\
			&\left\langle\left((u-w)+\frac{a_1w-a_2u}{t} \right)^q,v\right\rangle^{\beta_1-1},\\\\
			\frac{d}{dt}\psi_\gamma(t):=\frac{d}{dt}&\left\langle\left((u-w)+\frac{a_1w-a_2u}{t} \right)^\gamma,v\right\rangle^{\beta_2}=\\
			&-\gamma\beta_2\frac{1}{t^3}\left\langle\left((u-w)+\frac{a_1w-a_2u}{t} \right)^{\gamma-2}\left(t(u-w)+a_1w-a_2u \right)(a_1w-a_2u),v\right\rangle\times\\
			&\left\langle\left((u-w)+\frac{a_1w-a_2u}{t} \right)^\gamma,v\right\rangle^{\beta_2-1}.
		\end{align*}
		Assume that $(a_1w-a_2u)(x)\geq 0$ in $\Omega$. Note that by \eqref{talpha},  $t(u-w)+a_1w-a_2u =(a_1-a_2)(\alpha u+(1-\alpha)w)>0$. Hence
		\begin{equation}
				\frac{d}{dt}\psi_q(t)< 0, ~~\frac{d}{dt}\psi_q(t)< 0,~~\forall ~t \in [a_2,a_1],
		\end{equation}
		and thus
		$$
		(\psi_q(t)\psi_\gamma(t))'=\psi_q'(t)\psi_\gamma(t)+\psi_q(t)\psi_\gamma'(t)<0~~\Rightarrow~~\tilde{f}'(t)>0,~~~t \in [a_2,a_1].
		$$
		Similarly if $(a_1w-a_2u) \leq 0$ in $\Omega$, then $\tilde{f}'(t)<0$, $t \in [a_2,a_1]$.  
		Thus, $\tilde{f}(t)$ is a strong monotone  function on $[a_2,a_1]$ and we get the proof. 
	\end{proof}
	
	\begin{prop}\label{V11}
		Suppose that $a_1w(x)-a_2u(x) \neq 0$  in $\Omega$ and $(a_1w-a_2u)$ changes the sign in $\Omega$, then $f(\alpha) > \min\{f(0),f(1)\}$, $\forall \alpha \in [0,1]$.
	\end{prop}
	\begin{proof}
	  Suppose, contrary to our claim, that there exists $\alpha^* \in (0,1)$ such that
	\begin{equation}\label{ContrN1}
		f(\alpha^*)\leq \min\{f(0),f(1)\}\equiv \min\{\tilde{\lambda}(u,v)),\tilde{\lambda}(w,v))\}.
	\end{equation}
	Denote
	$
	\phi^*:=\alpha^* u+(1-\alpha^*)w
	$.
	It is easily seen that there exist $s_1>0$, $s_2>0$, $\tau^
	*\in [0,1]$ such that $\phi^*:=\tau^*s_1u+(1-\tau^*)s_2w$ and $(a_1 s_1 u-a_2 s_2 w) \geq 0,~\mbox{in}~\Omega$. 
Hence, by Proposition \ref{V1} and due to the homogeneity of $\min\{\tilde{\lambda}(s_1u,v)),\tilde{\lambda}(s_2w,v))\}$ we have
	$$
	f(\alpha^*) >\min\{\tilde{\lambda}(s_1u,v)),\tilde{\lambda}(s_2w,v))\}=\min\{\tilde{\lambda}(u,v)),\tilde{\lambda}(w,v))\}
	$$
	which contradicts \eqref{ContrN1}.
		\end{proof}

In the case $a_1-a_2= 0$, we have
	\begin{align*}
		f(\alpha):=&\frac{\left\langle -\Delta  w, v\right\rangle}{\left(\left\langle (\alpha u+(1-\alpha)w )^q,v\right\rangle\right)^{\beta_1}\left(\left\langle (\alpha u+(1-\alpha)w )^\gamma,v\right\rangle \right)^{\beta_2}}.
	\end{align*}

Similar to Propositions \ref{V1}, \ref{V11} we have 
	\begin{prop}\label{PV2}
	Suppose  $a_1-a_2= 0$ and $u-w \not\equiv 0$ in $\Omega$. Then $f(\alpha) > \min\{f(0),f(1)\}$, $\forall \alpha \in [0,1]$.
	\end{prop}
 Thus, Propositions \ref{V1}, \ref{V11}, \ref{PV2} yield  (i).

\medskip

	Let us  prove (ii). 	Observe
	$$
	\tilde{\lambda}( u,\alpha v+(1-\alpha)w)=\frac{\alpha\left\langle -\Delta  u, v\right\rangle+(1-\alpha)\left\langle -\Delta  u, w\right\rangle}{\left(\alpha\left\langle u^q,v\right\rangle +(1-\alpha)\left\langle u^q,w\right\rangle\right)^{\beta_1}\left(\alpha\left\langle u^\gamma,v\right\rangle +(1-\alpha)\left\langle u^\gamma,w\right\rangle\right)^{\beta_2}}
	$$
	Let us show
	\begin{equation}\label{GlIneq}
		\tilde{\lambda}( u,\alpha v+(1-\alpha)w)\leq \max\{\tilde{\lambda}(u,v)),\tilde{\lambda}(w,v))\}.
	\end{equation}
	Denote
	\begin{align*}
		s:=\alpha\left\langle u^q,v\right\rangle +(1-\alpha)\left\langle u^q,w\right\rangle,\\
		t:=\alpha\left\langle u^\gamma,v\right\rangle +(1-\alpha)\left\langle u^\gamma,w\right\rangle.	
	\end{align*}
	Assume first that $\left\langle u^q,v\right\rangle\left\langle u^\gamma,w\right\rangle-\left\langle u^q,w\right\rangle\left\langle u^\gamma,v\right\rangle\neq 0$.  Then there exist $x,y \in \mathbb{R}$ such that 
	\begin{equation}\label{std}
		\alpha\left\langle -\Delta  u, v\right\rangle+(1-\alpha)\left\langle -\Delta  u, w\right\rangle=sx+ty.
	\end{equation}
	Hence 
	$$
	\tilde{\lambda}( u,\alpha v+(1-\alpha)w)=\frac{sx+ty}{s^{\beta_1}t^{\beta_2}}=x (\frac{s}{t})^{\beta_2}+y (\frac{t}{s})^{\beta_1}=x r^{\beta_2}+y r^{-\beta_1}=:g(r).
	$$
	where $r=\frac{s}{t}$.
	Calculate
	\begin{align}
		g'(r)=&\beta_2 x r^{\beta_2-1}-\beta_1 y r^{-\beta_1-1}=r^{-\beta_1-1}(\beta_2 xr-\beta_1 y),\label{eq:g}\\
		g''(r)=&\beta_2 (\beta_2-1) x r^{\beta_2-2}+\beta_1 (\beta_1+1) y r^{-\beta_1-2}=
		r^{-\beta_1-2}\left(\beta_2 (\beta_2-1) xr+\beta_1 (\beta_1+1) y \right). \nonumber
	\end{align}
	Hence $g'(\bar{r})=0$ if and only if
	 $xy>0$. Since $ \bar{r}=\frac{\beta_1 y}{\beta_2 x}$,
	$$
	g''(\bar{r})=\bar{r}^{(-\beta_1-2)}\left( \beta_1(\beta_2-1)y +\beta_1 (\beta_1+1) y \right)=\beta_1\bar{r}^{(-\beta_1-2)}y.
	$$
	If $y<0$, $x<0$, then $sx+ty<0$ which is impossible. Thus  $y>0$, and 
	$	g''(\bar{r})>0$, which implies \eqref{GlIneq}. 
	If $xy\leq 0$,  \eqref{eq:g} implies that $g(r)$ is monotone function and thus, 
	\eqref{GlIneq} holds true.
	
In the  case $\left\langle u^q,v\right\rangle\left\langle u^\gamma,w\right\rangle-\left\langle u^q,w\right\rangle\left\langle u^\gamma,v\right\rangle= 0$, we have $s=\kappa t$ for some $\kappa \in \mathbb{R}_+\setminus 0$, and consequently,
$$
f(\alpha):=\frac{\alpha\left\langle -\Delta  u, v\right\rangle+(1-\alpha)\left\langle -\Delta  w, v\right\rangle}{\kappa^{\beta_1}(\alpha\left\langle u^\gamma,v\right\rangle +(1-\alpha)\left\langle u^\gamma,w\right\rangle) }, ~~\alpha \in [0,1],~u, v \in S.
$$
Note that this is a fractional linear function on $\alpha$. Since $\alpha\left\langle u^\gamma,v\right\rangle +(1-\alpha)\left\langle u^\gamma,w\right\rangle>0$ for $\alpha \in [0,1]$, this implies that $f(\alpha)$ is a monotone function and thus, \eqref{GlIneq} holds true.


\begin{thebibliography}{25}





\bibitem{ABC}  A. Ambrosetti, H. Brezis,   G. Cerami,   Combined effects of concave and convex nonlinearities in some elliptic problems, J. Funct. Anal. 122(2) (1994) 519--543.


\bibitem{Bagirov} A. Bagirov,  N. Karmitsa, \& M. M. M\"akel\"a,   Introduction
to Nonsmooth Optimization: Theory, Practice
and Software (Springer), 2014.












\bibitem{Cazenave} T. Cazenave,  M.  Escobedo,  M.A. Pozio,  Some stability properties for minimal solutions of $-\Delta u= \lambda g (u)$, Portugaliae Mathematica, 59, (2002) 373--391.


\bibitem{Gidas} B. Gidas, W.-M. Ni, L. Nirenberg, Symmetry and related properties via the
maximum principle, Commun. Math. Phys. 68,  (1979) 209-243.

\bibitem{CranRibin}  M. G. Crandall, P. H. Rabinowitz,  Bifurcation, perturbation of simple eigenvalues, and linearized stability. Arch. Rat. Mech. Anal. 52(2), (1973) 161--180.



	
%
%
\bibitem{Dieudonn}J. Dieudonn\'e, Foundations of Modern Analysis. Academic Press, New York,
1964.



\bibitem{IlyasJDE24}  Y. Il'yasov, A finding of the maximal saddle-node bifurcation for systems of differential equations. Journal of Differential Equations 378, (2024) 610-625.


















\bibitem{Giltrud} D. Gilbarg, N.S. Trudinger,  Elliptic partial differential equations of second order. Berlin: Springer, 1977.









\bibitem{IlyasChaos} 
Y. Ilyasov,  Finding Saddle-Node Bifurcations via a Nonlinear Generalized Collatz–Wielandt Formula. International Journal of Bifurcation and Chaos. 2021 31(01):2150008.












\bibitem{IlyasFunc}
Y.~S. Il’yasov, Bifurcation calculus by the extended functional method,   Functional Analysis and Its Applications, 41(1)1, (2007) 18--30.


\bibitem{IvanIlya} 
Y.S. Il'yasov,  A.A. Ivanov, Finding bifurcations for solutions of nonlinear equations by quadratic programming methods,  { Comp. Math. and Math. Phys.}, 53(3), (2013) 350--364.

\bibitem{IlIvan1}
Y.S. Il'yasov,   A.A. Ivanov, Computation of maximal turning points to nonlinear   equations by nonsmooth optimization, Optim. Meth. and   Softw., 31 (1), ( 2016) 1--23.











\bibitem{keller1}
H.B. Keller,  Numerical solution of bifurcation and nonlinear eigenvalue problems,
in Application of bifurcation theory, ed. P. Rabinowitz ( Academic Press, New York), (1977) 359--384.




\bibitem{kielh} H. Kielh\"{o}fer,  Bifurcation theory: An introduction with applications to PDEs. Vol. 156. Springer Science \& Business Media, 2006.



\bibitem{Korman} P. Korman,  On uniqueness of positive solutions for a class of semilinear equations. Discrete and Continuous Dynamical Systems, 8(4), (2002) 865-872.

\bibitem{Korman2} P. Korman, Global solution curves for semilinear elliptic equations. World Scientific, 2012.

\bibitem{Kuznet} Y.A. Kuznetsov, Elements of Applied Bifurcation Theory, vol.112, Springer Science \& Business Media, 2013.



\bibitem{Ouyang} T. Ouyang,  and J. Shi, Exact multiplicity of positive solutions for a class of semilinear problem, II. Journal of Differential Equations, 158(1), (1999) 94-151.


\bibitem{Tang} M.Tang, Exact multiplicity for semilinear elliptic Dirichlet problems involving concave and convex nonlinearities. Proceedings of the Royal Society of Edinburgh Section A: Mathematics, 133(3), ( 2003) 705-717.

\bibitem{Salazar} P. D. P. Salazar, Y. Il'yasov, L. F. C. Alberto,  E. C. M. Costa \& M. B. Salles,  Saddle-node bifurcations of power systems in the context of variational theory and nonsmooth optimization, IEEE Access 8 (2020): 110986--110993.


\bibitem{Sion}  M. Sion,  On general minimax theorems. (1958): 171-176.

	
\end{thebibliography}
\end{document}